\documentclass[10pt,a4paper]{article}
\usepackage{multicol}

\begin{document}

\setcounter {equation}{0}

{\bf GEOMETRIC CHARACTERIZATION OF STRICTLY POSITIVE REAL REGIONS
AND ITS APPLICATIONS}

{\vspace{0.34cm}  
Long Wang and Wensheng Yu

{\vspace{0.37cm} Center for Systems and Control, Peking
University, Beijing, China

\vspace{0.87cm}

\begin{multicols}{2}
{\bf Abstract:} { Strict positive realness (SPR) is an important
concept in absolute stability theory, adaptive control, system
identification, etc. This paper characterizes the strictly
positive real (SPR) regions in coefficient space and presents a
robust design method for SPR transfer functions. We first
introduce the concepts of SPR regions and weak SPR regions and
show that the SPR region associated
 with a fixed polynomial is unbounded, whereas the weak SPR region is bounded. We then
prove that the intersection of several weak SPR regions associated
with different polynomials can not be a single point. Furthermore,
we show how to construct a point in the SPR region from a point in
the weak SPR region. Based on these theoretical development, we
propose an algorithm for robust design of SPR transfer functions.
This algorithm works well for both low order and high order
polynomial
 families. Illustrative examples are provided to show the effectiveness of this algorithm.
}


\vskip 0.5cm
\section{\large   Introduction}

The notion of strict positive realness (SPR) of transfer functions
plays an important role in absolute stability theory, adaptive
control and system identification[1-23]. Motivated by Kharitonov's
seminal theorem on the robust stability for a family of
polynomials, a number of recent papers has concentrated  on the
strict positive realness for a family of transfer functions. In
the spirit of Kharitonov, the robust SPR analysis and design
problems were first formulated by Dasgupta and Bhagwat[8]. They
showed that every transfer function in an interval transfer
function family is strictly positive real if and only if sixteen
prescibed vertex transfer functions in this family are strictly
positive real. The sixteen critical vertex transfer functions can
be constructed explicitly using Kharitonov's four vertex
polynomials. This result was subsequently improved by Chapellat
and Bhattacharyya, Wang and Huang, where only eight out of the
sixteen critical vertex transfer functions need to be
checked[11-12]. For a family of transfer functions with affine
linearly correlated perturbations, or more generally,
multilinearly correlated perturbations, Dasgupta, Anderson et al.
showed that it suffices to check all vertices in order to ensure
the strict positive realness of the entire family[14]. By resort
to the concept of positive polynomial pairs and root interlacing
properties, Hollot and Huang solved the robust SPR design problem
for low order and structured families[9-10]. Anderson et al.
considered the general robust SPR design problem, and by using the
Hilbert transform, provided a constructive method[14]. Betser and
Zeheb made some further improvements[15].

This paper characterizes SPR regions in coefficient space and
presents a robust design method for SPR transfer functions. We
first introduce the concepts of SPR regions and weak SPR regions
and give a complete characterization of them. We show that the SPR
region associated
 with a fixed polynomial is unbounded, whereas the weak SPR region is bounded. We then
prove that the intersection of several weak SPR regions associated
with different polynomials can not be a single point. Furthermore,
we show how to construct a point in the SPR region from a point in
the weak SPR region. Based on these theoretical development, we
propose an algorithm for robust design of  SPR transfer functions.
This algorithm works well for both low order and high order
polynomial
 families. Illustrative examples are provided to show the effectiveness of this algorithm.

\vskip 0.7cm
\section{\large   Preliminaries}

Denote $ P^n$ as the $n$-th order real polynomial family, $ R^n$
as the $n$ dimensional real field, and $H^n\subset P^n$ as the set
of all $n$-th order Hurwitz  stable polynomials.

In the following definitions,$ b(\cdot )\in P^m, a(\cdot )\in
P^n,$ and $p(s)=b(s)/a(s)$ is a rational function.

{\bf Definition 1}\ \ $p(s)$ is said to be strictly positive
real(SPR),denoted as $p(s)\in$SPR,if $b(s)\in P^n,a(s)\in H^n,$
and Re$[p(j\omega )]>0,$ $\forall \omega \in R.$

{\bf Definition 2}\ \ $p(s)$ is said to be weak SPR (WSPR),denoted
as $p(s)\in$WSPR,if $b(s)\in P^{n-1},a(s)\in H^n,$ and
Re$[p(j\omega )]>0,$ $\forall \omega \in R.$

{\bf Definition 3}\ \ Given $a(s)\in H^n,$ the set of the
coefficients (in $R^{n+1}$) of all the $b(s)$'s in $P^n$ such that
$p(s) :=\displaystyle\frac {b(s)}{a(s)}\in$SPR is said to be the
SPR region associated with $a(s)$,denoted as $\Omega _{a}.$

{\bf Definition 4}\ \ Given $a(s)\in H^n,$ the set of the
coefficients (in $R^{n}$) of all the $b(s)$'s in $P^{n-1}$ such
that $p(s) :=\displaystyle\frac {b(s)}{a(s)}\in$WSPR is said to be
the WSPR region associated with $a(s)$,denoted as $\Omega
^{W}_{a}.$

For notational convenience, $\Omega _{a}(\Omega ^{W}_{a})$
sometimes also stands for the set of all the polynomials $b(s)$ in
$P^n(P^{n-1})$ such that $p(s) :=\displaystyle\frac
{b(s)}{a(s)}\in$SPR(WSPR).

From the definitions above, it is easy to get the following
properties:

{\bf Property 1}$^{[1,4,9,14,19]}$\ \ If $p(s)\in$SPR(WSPR), then
$|arg(b(j\omega ))- arg(a(j\omega ))|<\displaystyle\frac {\pi }{2}
$ ,  $\forall \omega \in R,$ where $arg(\cdot)$ stands for the
argument of the complex number,and the difference of two arguments
can differ by an integer number of $2\pi$.

{\bf Property 2}$^{[9, 10]}$\ \ Given $a(s)\in H^n,$ $\Omega _{a}$
is a non-empty, open, convex cone in $R^{n+1}$.

{\bf Property 3}$^{[10, 11]}$\ \ Given $a(s)\in H^n,$ we have
$\Omega _{a}\subset H^n,
\Omega _{a}^W\subset H^{n-1}.$%

The problem we are interested in is: Given a family of Hurwitz
stable polynomials, how can find a fixed polynomial such that
their ratios will be SPR-invariant? In what follows, we will first
give some characterization of WSPR regions, and then propose an
efficient design procedure for this problem.

\vskip 0.7cm
\section{\large   Geometric Characterization of SPR Regions}

By definition, an SPR (WSPR) transfer function times a positive
integer is still SPR (WSPR). Thus, without loss of generality, let

\begin{equation}
a(s)=s^n+a_1s^{n-1}+\cdots +a_n\in H^n
\end{equation}

Denote as $\Omega _{1a}$ the set of the coefficients of all the
$b(s)=s^n+x_1s^{n-1}+\cdots +x_n\in P^n,$ i.e.,
$(x_1,x_2,\cdots,x_n)$ in $R^{n}$ , such that $p(s)
=\displaystyle\frac {b(s)}{a(s)}\in$SPR; and denote as $\Omega
^W_{1a}$  the set of the coefficients of all the
$b(s)=s^{n-1}+x_1s^{n-2}+\cdots +x_{n-1}\in P^{n-1},$ i.e.,
 $(x_1,x_2,\cdots,x_{n-1})$ in $R^{n-1}$ , such that $p(s)
=\displaystyle\frac {b(s)}{a(s)}\in$WSPR. Obviously, we have

\begin{eqnarray}
\{1\}\times \Omega _{1a}=\bigg\{(1, \displaystyle\frac
{b_1}{b_0}, \displaystyle\frac {b_2}{b_0},\cdots,
 \displaystyle\frac {b_n}{b_0})\bigg | \hspace{2.5cm}\nonumber\\
\forall (b_0,b_1,b_2,\cdots,b_n) \in \Omega _{a}\bigg\} \\
\{1\}\times \Omega ^W_{1a}=\bigg \{(1, \displaystyle\frac
{b_1}{b_0}, \displaystyle\frac {b_2}{b_0},\cdots,
 \displaystyle\frac {b_{n-1}}{b_0}) \bigg |\hspace{2cm}\nonumber\\
\forall (b_0,b_1,b_2,\cdots,b_{n-1}) \in \Omega ^W_{a} \bigg \}
\end{eqnarray}

For notational convenience, $\Omega _{1a}(\Omega ^{W}_{1a})$
sometimes also stands for the corresponding polynomial set.

As we know$^{[9, 10]}$, $\Omega _{a}$ is a non-empty, open, convex
cone in $R^{n+1}$. Thus, $\Omega _{a}$ is an unbounded set in
$R^{n+1}$. In what follows, we will show that $\Omega _{1a}$ is
also an unbounded set in $R^{n}$.

{\bf Theorem 1}\ \ Given $a(s)\in H^n,$ $\Omega _{1a}$ is a
non-empty, open, unbounded convex set in $R^{n}$.

{\bf Proof}\ \ Obviously, we have $a(s)\in \Omega _a.$ If the
leading coefficient of $a(s)$ is $a_0,$ then

\begin{center}
 $(\displaystyle\frac{a_1}{a_0},
\displaystyle \frac{a_2}{a_0},\cdots , \displaystyle
\frac{a_n}{a_0})\in \Omega _{1a}.$
\end{center}

\noindent Thus, $\Omega _{1a}$ is not empty.

Moreover, $(1,x_1,x_2,\cdots ,x_n)\in \Omega _a, $ $\forall
(x_1,x_2,\cdots ,x_n)\in \Omega _{1a}.$ By Property 2, $\Omega _a$
is open. Thus, these exists $\delta >0,$ such that, when
$\sqrt{(1-y_0)^2+(x_1-y_1)^2+\cdots +(x_n-y_n)^2}<\delta $, we
have $(y_0,y_1,y_2,\cdots ,y_n)\in \Omega _a.$ For this $\delta ,$
if $(z_1,z_2,\cdots ,z_n)\in R^n$ satisfies

\begin{center}
$\sqrt{(x_1-z_1)^2+\cdots +(x_n-z_n)^2}<\delta $.
\end{center}

\noindent Then, we have
$$\sqrt{(1-1)^2+(x_1-z_1)^2+\cdots+(x_n-z_n)^2}<\delta .$$
 Hence $(1,z_1,z_2,\cdots ,z_n)\in \Omega
_a.$ Thus, we have $(z_1,z_2,\cdots ,z_n)\in \Omega _{1a}.$
Namely, $\Omega _{1a}$ is an open set.

By definition, $\Omega _{1a}$ is convex.

In what follows, we will prove that $\Omega _{1a}$ is unbounded.
For this purpose, we first introduce some notations, which are
needed in other proofs as well. Let
\begin{equation}
a(s)=s^n+a_1s^{n-1}+\cdots +a_n\in H^{n}
\end{equation}

\begin{equation}
b(s)=x_0s^{n}+x_1s^{n-1}+\cdots +x_{n}\in P^n\cup P^{n-1}
\end{equation}
Then $\forall \omega \in R,$ we have

\begin{center}
$\begin{array}{ll} \mbox{Re}[\displaystyle\frac {b(j\omega
)}{a(j\omega )}] &=\displaystyle\frac {1}{|a(j\omega
)|^2}\mbox{Re}[b(j\omega )a(-j\omega )]\\
&=\displaystyle\frac {1}{|a(j\omega
)|^2}\sum\limits_{l=0}^{n}(\sum\limits_{k=0}^{n}a_{k}
x_{2l-k}(-1)^{l+k})\omega
^{2(n-l)}\\
&=\displaystyle \frac {1}{|a(j\omega )|^2}\sum\limits_{l=0}^{n}
c_l\omega ^{2(n-l)}\\
\end{array}$
\end{center}
where $c_l:=\sum\limits_{k=0}^{n}a_{k}x_{2l-k}(-1)^{l+k},$ where
$a_0=1,$  and let $a_i=0, x_i=0,$  when $i<0$ or $i>n,$
$l=0,1,\cdots,n.$

Introducing the matrices
$$H_a:=\left[
\setlength{\arraycolsep}{0.035cm}
\begin{array}{ccccccc}
a_1 & 1 & 0 & 0 & 0 & \cdots  & 0 \\
a_3 & a_2 & a_1 & 1 & 0 & \cdots  & 0 \\
a_5 & a_4 & a_3 & a_2 & a_1 & \cdots  & 0 \\
\vdots  & \vdots  & \vdots  & \vdots  & \vdots  & \ddots  & \vdots  \\
a_{2n-1} & a_{2n-2} & a_{2n-3} & a_{2n-4} & a_{2n-5} & \cdots  &
a_n
\end{array}
\right] ,$$
$$E_n:=\left[
\begin{array}{ccccc}
1 &  &  &  &  \\
& -1 &  &  &  \\
&  & 1 &  &  \\
&  &  & -1 &  \\
&  &  &  & \ddots
\end{array}
\right] ,$$

$$A:=\left[
\begin{array}{cccc}
1 & 0 & \cdots  & 0 \\
-a_2 &  &  &  \\
a_4 &  & E_nH_aE_n &  \\
\vdots  &  &  &
\end{array}
\right] ,$$
$$b:=\left[
\begin{array}{c}
x_0 \\
x_1 \\
\vdots  \\
x_{n-1}
\end{array}
\right] ,$$
$$c:=\left[
\begin{array}{c}
c_0 \\
c_1 \\
\vdots  \\
c_n
\end{array}
\right] $$
where $a_i=0$ when $i>n.$ Then, it is easy to verify that

\begin{equation}
c=Ab
\end{equation}

Since $a(s)\in H^n,$ we know that $E_nH_aE_n$ is invertible. For
any $d=[d_1,d_2,\cdots ,d_n]^T\in R^n$ such that all elements of
$d$ are positive. Denote $\bar a=[-a_2,a_4,-a_6,a_8,\cdots
,(-1)^na_{2n}]^T,$ where $a_i=0$ when $i>n$. Let $\bar
b=(E_nH_aE_n)^{-1}(d-\bar a):=[b_1,b_2,\cdots, b_n]^T,$ Then
obviously, we have $[1,d_1,d_2,\cdots ,d_n]^T=A[1,b_1,b_2,\cdots
,b_n]^T,$ and $b_n=\displaystyle\frac{d_n}{a_n}.$ By $c=Ab$, we
know that $(b_1,b_2,\cdots, b_n)\in \Omega _{1a}.$

On the other hand, due to the arbitrariness of $d$, $d_n$ can be
taken arbitrarily large. Therefore, $b_n$ can also be arbitrarily
large. Namely, $\Omega _{1a}$ is unbounded. This completes the
proof.

{\bf Remark 1} Given two stable polynomials $a_1 (s)$ and $a_2
(s)$ , existence of a polynomial $b (s)$ such that $\frac{b (s)}{
a_1 (s) }$ and
 $\frac{b (s)}{ a_2 (s) }$ are both SPR is tantamount to non-emptiness of the
 intersection of the two SPR regions associated with  $a_1 (s)$ and $a_2 (s)$ .
Since $\Omega _{a}$ and $\Omega _{1a}$ are both unbounded sets.
When dealing with robust SPR design problem, we must find the
intersection of several unbounded sets (i.e., SPR regions), which
is  intractable. This is the reason that we introduce the concept
of WSPR regions, which are bounded as shown below.

{\bf Theorem 2}\ \ Given $a(s)\in H^n,$ $\Omega ^W_{1a}$ is a
non-empty, bounded convex set in $R^{n-1}$.

{\bf Proof}\ \

$H_a,E_n$ and $A$ were defined in the proof of Theorem 1.

Denote $B$ as the $(n-1)\times (n-1)$ matrix formed by the first
$n-1$ row and last $n-1$ column of the matrix $E_n H_a E_n$.
Obviously, $B$ is also invertible.

Denote $\bar a:=[a_1,-a_3,a_5,-a_7,\cdots
,(-1)^{n-1}a_{2(n-1)+1}]^T$ ( $a_i=0$ when $i>n$).

Let $\bar b:=-B^{-1}\bar a=[b_1,b_2,\cdots ,b_{n-1}]^T.$ Since
$a(s)\in H^n,$ it is easy to verify that $b_{n-1}>0.$ Denote
$b=[0,1,b_1,b_2,\cdots ,b_{n-1}]^T.$ Let $c_0=c_1=\cdots
=c_{n-1}=0,$ $c_n=a_nb_{n-1},$ in $c:=[c_0,c_1,\cdots ,c_n]$.
Then, it is easy to verify that $c=Ab$ is true. Thus, we have
$(b_1,b_2,\cdots ,b_{n-1})\in \Omega _{1a}^W,$ Namely, $\Omega
_{1a}^W$ is not empty.

By definition, $\Omega _{1a}^W$ is convex.

We now prove that $\Omega _{1a}^W$ is bounded.

For any $(x_1,x_2,\cdots ,x_{n-1})\in \Omega _{1a}^W,$ we have
$$\displaystyle\frac{s^{n-1}+x_1s^{n-2}+ \cdots +x_{n-1}}
{s^n+a_1s^{n-1}+\cdots +a_n}\in
 WSPR,$$
By Property 3, $s^{n-1}+x_1s^{n-2}+\cdots +x_{n-1}\in H^{n-1}.$
Moreover, $\forall \omega \in R,$ we have

\begin{equation}
\mbox{Re}(\displaystyle \frac{s^{n-1}+x_1s^{n-2}+\cdots +x_{n-1}}
{s^n+a_1s^{n-1}+\cdots +a_n}%
|_{s=j\omega })>0.
\end{equation}

\noindent Thus

\begin{equation}
\mbox{Re}(\displaystyle \frac{s^n+a_1s^{n-1}+\cdots +a_n}
{s^{n-1}+x_1s^{n-2}+\cdots +x_{n-1}}%
|_{s=j\omega })>0,\ \ \ \forall \omega \in R.
\end{equation}

\noindent Obviously
{\scriptsize{
\begin{eqnarray*}
\frac{s^n+a_1s^{n-1}+\cdots +a_n}{s^{n-1}+x_1s^{n-2}+\cdots +x_{n-1}}=s+\hspace{4cm}\nonumber\\
\frac{(a_1-x_1)s^{n-1}+(a_2-x_2)s^{n-2}+\cdots +(a_{n-1}-x_{n-1})s+a_n}{s^{n-1}+x_1s^{n-2}+\cdots +x_{n-1}}
\end{eqnarray*}
}}
\noindent Therefore

{\scriptsize{
\begin{eqnarray*}
Re(\frac{(a_1-x_1)s^{n-1}+(a_2-x_2)s^{n-2}+\cdots +(a_{n-1}-x_{n-1})s+a_n}{%
s^{n-1}+x_1s^{n-2}+\cdots +x_{n-1}}|_{s=j\omega })\\
=Re(\displaystyle \frac{s^n+a_1s^{n-1}+\cdots +a_n}
{s^{n-1}+x_1s^{n-2}+\cdots +x_{n-1}}%
|_{s=j\omega })-Re(j\omega )\hspace{2.5cm}\\
=\mbox{Re}(\displaystyle \frac{s^n+a_1s^{n-1}+\cdots +a_n}
{s^{n-1}+x_1s^{n-2}+\cdots +x_{n-1}}%
|_{s=j\omega })>0,\ \ \ \forall \omega \in R.\hspace{1.8cm}\\
\end{eqnarray*}
}}
\noindent It is easy to see that

{\scriptsize{
$$
\frac{(a_1-x_1)s^{n-1}+(a_2-x_2)s^{n-2}+\cdots +(a_{n-1}-x_{n-1})s+a_n}{%
s^{n-1}+x_1s^{n-2}+\cdots +x_{n-1}}$$
$$ \hspace{5cm}\in \{\mbox{SPR}\}\cup
\{\mbox{WSPR}\}.$$
}}
\noindent Again, by Property 3, we have
\begin{eqnarray}
(a_1-x_1)s^{n-1}+(a_2-x_2)s^{n-2}+\cdots \hspace{1.5cm}\nonumber\\
+(a_{n-1}-x_{n-1})s+a_n\in H^{n-1}\cup H^{n-2}
\end{eqnarray}

\noindent Hence

\begin{equation}
0<x_1\leq a_1,0<x_2<a_2,\cdots ,0<x_{n-1}<a_{n-1}
\end{equation}

\noindent Namely

$$
\Omega _{1a}^W\subset \{(x_1,x_2,\cdots ,x_{n-1})|\alpha(s):=\hspace{3cm}$$
$$\sum_{i=1}^n(a_i-x_i)s^{n-i}
\in H^{n-1}\cup H^{n-2}, \mbox{where}\,\,\,\,\,  x_n=0\}$$
$$\subset \{(x_1,x_2,\cdots ,x_{n-1})|0<x_1\leq
a_1,0<x_2<a_2,\cdots,$$
$$\hspace{5cm}0<x_{n-1}<a_{n-1}\}.$$

\noindent Thus, $\Omega _{1a}^W$ is bounded. This completes the
proof.

{\bf Remark 2} It should be pointed out that $\Omega ^W_{1a}$ is
not an open set in $R^{n-1}$. In fact, from the proof of Theorem
2, we know that $\Omega ^W_{1a}$ is tangent to the hyperplane
$x_1=a_1$ in $R^{n-1}$. And there exist some points of $\Omega
^W_{1a}$ in this hyperplane. Thus, $\Omega ^W_{1a}$ can not be an
open set. Obviously,$\Omega ^W_{a}$ is a non-empty, convex cone in
$R^{n-1}$. Thus, $\Omega ^W_{a}$ is also unbounded. One may be
tempted to believe that $\Omega ^W_{a}$ is not an open set either.

Though $\Omega ^W_{1a}$ is not an open set. The following theorem
guarantees such a fact: when the intersection of two or more WSPR
regions is not empty, then the intersection must be a region, not
a single point. This means that Ackermann's counterexample (that
the unstable region is an isolated point so that gridding the
parameter space can lead to erroneous conclusions no matter how
dense the gridding is[19]) does not happen in this case.

{\bf Theorem 3}\ \ Given $a(s)\in H^n,$ if
$(x_1,x_2,\cdots,x_{n-1})\in \Omega ^W_{1a},$ then for
sufficiently small $\varepsilon>0,$ we have
$(x_1-\varepsilon,x_2-\varepsilon, \cdots,x_{n-1}-\varepsilon)\in
\Omega ^W_{1a}.$

{\bf Proof}\ \ $\forall (x_1,x_2,\cdots ,x_{n-1})\in \Omega
_{1a}^W,$ and $\forall \omega \in R,$ we have

\begin{equation}
\mbox{Re} (\displaystyle
\frac{s^{n-1}+x_1s^{n-2}+\cdots +x_{n-1}}{s^n+a_1s^{n-1}+\cdots +a_n}%
|_{s=j\omega })>0.
\end{equation}

\noindent $\forall \varepsilon >0,$ since
$$
Re(\frac{s^{n-1}+(x_1-\varepsilon )s^{n-2}+\cdots +(x_{n-1}-\varepsilon )}{
s^n+a_1s^{n-1}+\cdots +a_n}|_{s=j\omega })$$
$$=Re(\frac{s^{n-1}+x_1s^{n-2}+\cdots +x_{n-1}}{s^n+a_1s^{n-1}+\cdots +a_n}|_{s=j\omega }) \hspace{2.8cm}$$
$$\hspace{1cm}+Re( \frac{(-\varepsilon )(s^{n-2}+s^{n-3}+\cdots +1)}{
s^n+a_1s^{n-1}+\cdots +a_n}|_{s=j\omega })$$
$$=Re(\frac{s^{n-1}+x_1s^{n-2}+\cdots +x_{n-1}}{s^n+a_1s^{n-1}+\cdots +a_n}|_{s=j\omega })\hspace{2.8cm}$$
$$+ \frac{(-\varepsilon )}{|a(j\omega
)|^2}(-\omega
^{2(n-1)}+\tilde c(\omega )),\\
$$

where $\tilde c(\omega )$ is a real polynomial of order less or
equal to $2(n-2)$. Thus, when $\mid \omega \mid$ is sufficiently
large, the sign of $(-\varepsilon ) (-\omega ^{2(n-1)}+\tilde
c(\omega ))$ will be positive. Namely, there exists $ \omega _1
>0$ such that, for all $ \mid \omega \mid \geq \omega _1,$

\begin{equation}
\mbox{Re}(\displaystyle\frac
{s^{n-1}+(x_1-\varepsilon )s^{n-2}+\cdots +(x_{n-1}-\varepsilon )}{%
s^n+a_1s^{n-1}+\cdots +a_n}|_{s=j\omega })
>0.
\end{equation}

Denote

\begin{equation}
 M_1=\inf_{\mid \omega \mid \leq \omega _1}\mbox{Re}( \displaystyle\frac
{ s^{n-1}+x_1s^{n-2}+\cdots +x_{n-1}} { s^n+a_1s^{n-1}+\cdots +a_n
}|_{s=j\omega }),
\end{equation}

\begin{equation}
 N_1=\sup_{\mid \omega \mid \leq \omega _1}\displaystyle {\mid
\mbox{Re}(\displaystyle\frac {1}{|a(j\omega )|^2}(\omega
^{2(n-1)}-\tilde c(\omega )) \mid }
\end{equation}

\noindent Then $M_1>0 $ and $N_1>0.$ Choosing
$0<\varepsilon<\displaystyle\frac {M_1}{N_1},$ by simple
computation, we have
{\scriptsize{
\begin{equation}
 \mbox{Re}
(\displaystyle\frac
{s^{n-1}+(x_1-\varepsilon )s^{n-2}+\cdots +(x_{n-1}-\varepsilon )}{%
s^n+a_1s^{n-1}+\cdots +a_n}|_{s=j\omega })
>0, \forall \omega \in R
\end{equation}
}}
\noindent Therefore,
$$\displaystyle\frac
{s^{n-1}+(x_1-\varepsilon )s^{n-2}+\cdots +(x_{n-1}-\varepsilon )}{%
s^n+a_1s^{n-1}+\cdots +a_n}\in WSPR,$$
 namely
$(x_1-\varepsilon,x_2-\varepsilon, \cdots,x_{n-1}-\varepsilon)\in
\Omega ^W_{1a}.$ This completes the proof.

The following theorem shows the relationship between $\Omega
^W_{1a}$ and $\Omega _{a}$, and plays an important role in robust
SPR design.

{\bf Theorem 4}\ \ Given $a(s)\in H^n,$ if
$(x_1,x_2,\cdots,x_{n-1})\in \Omega ^W_{1a},$ then $\forall
(1,\alpha_1,\alpha_2,\cdots, \alpha_n)\in R^{n+1},$ we can take
sufficiently small $\varepsilon>0$ such that $(0,1,x_1,x_2,
\cdots,x_{n-1})+\varepsilon (1,\alpha_1,\alpha_2,\cdots,
\alpha_n)\in \Omega _{a}.$

{\bf Proof}\ \ Denote

\begin{equation}
b(s)=s^{n-1}+x_1s^{n-2}+\cdots +x_{n-1},
\end{equation}

\begin{equation}
\alpha(s)=s^{n}+\alpha_1s^{n-1}+\cdots +\alpha_{n},
\end{equation}

\begin{equation}
{\tilde {b}}(s)=b(s)+\varepsilon \alpha(s).
\end{equation}

\noindent Since $(x_1,x_2,\cdots,x_{n-1})\in \Omega ^W_{1a},$ we
have

\begin{equation}
 \mbox{Re}
(\displaystyle\frac {b(j\omega )}{a(j\omega )})>0,  \forall \omega
\in R.
\end{equation}

We only need to show that, for sufficiently small $\varepsilon >0,
$

\begin{equation}
\mbox{Re} (\displaystyle\frac {\tilde b (j\omega )}{a(j\omega
)})>0,  \forall \omega \in R
\end{equation}

\noindent Obviously, $ \tilde {b}(s)$ and $ a(s)$ have same order
$n$. Thus, there exists $ \omega _2 >0$ such that, for all $ \mid
\omega \mid \geq \omega _2,$ we have Re$ (\displaystyle\frac
{\tilde {b}(j\omega )}{a(j\omega )})>0.$

Denote

\begin{equation}
 M_2=\inf_{\mid \omega \mid \leq \omega _2}\mbox{Re} ( \displaystyle\frac
{b(j\omega )}{a(j\omega )}),  \ \  \ \
 N_2=\sup_{\mid \omega \mid \leq \omega _2}\displaystyle {\mid
\mbox{Re}(\displaystyle\frac {\alpha(j\omega )}{a(j\omega )}) \mid
}
\end{equation}

\noindent Then $M_2>0 $ and $N_2>0.$ Choosing $0<\varepsilon
<\displaystyle\frac {M_2}{N_2},$ by simple computation, we have

\begin{equation}
 \mbox{Re}
(\displaystyle\frac {\tilde {b}(j\omega )}{a(j\omega )})>0,
\forall \omega \in R
\end{equation}

\noindent This completes the proof.

\vskip 0.7cm
\section{\large    Applications in Robust Design of SPR Transfer Functions}

Generally speaking, the design problem is more difficult than
analysis problem, since it is usually constructive, i.e., it not
only shows the existence of the solution, but also provides a
constructive procedure to find it. In this section, we will
propose an algorithm for robust design of  SPR transfer functions.
This algorithm works well for both low order and high order
polynomial
 families. Illustrative examples are provided to show the effectiveness of this algorithm.

Suppose

\begin{equation}
F=\{a_i(s)=s^n+\sum_{l=1}^na_l^{(i)}s^{n-l},i=1,2,\cdots,m.\}
\end{equation}
How can we find a polynomial $b(s),$ such that
$p_i(s):=\displaystyle\frac
{b(s)}{a_i(s)}\in$SPR,$i=1,2,\cdots,m$?

As observed earlier, existence of such a polynomial $b (s)$ boils
down to the condition that the intersection of the SPR regions
associated with $a_i (s)$  is not empty. From the results in the
previous section, we know that SPR regions are unbounded, whereas
WSPR regions are bounded. Thus, by a computational
 consideration, we first consider the
intersection of WSPR regions, and then construct a polynomial  $b
(s)$ by using the technique presented in the previous section.

Since SPR(WSPR) transfer functions with fixed numerator (or
denominator) enjoy convexity property, namely, if there exists a
polynomial $c(s)$, such that $\displaystyle\frac {c(s)}{a(s)}$ and
$\displaystyle\frac {c(s)}{b(s)}$ are both SPR(WSPR),then, it is
easy to verify that, for any $\alpha \geq 0,\beta \geq 0 $ and
$(\alpha ,\beta )\neq (0,0),$ we have $\displaystyle\frac
{c(s)}{\alpha a(s)+\beta b(s)}\in $SPR(WSPR). Therefore, the
assumptions made on $F$ do not lose any generality. Actually, the
method proposed in our paper also applies to convex combination of
polynomials, interval polynomials, and more generally, polytopic
polynomials and multilinearly perturbed polynomials[1,4,14].

By the results presented in the previous section, we propose the
following design procedure:

{\bf Step 1.} Test the robust stability of the convex hull of $F$,
i.e.,  $\overline{F}$. If $\overline{F}$ is robustly stable, then
go to Step 2; otherwise, print "there does not exist such a $b(s)$
"; (by Definitions 1 and 2)

{\bf Step 2.} Let $\alpha_l=\min
\{a_l^{(i)},i=1,2,\cdots,m\},l=1,2,\cdots,n-1.$ Gridding the
hyperrectangle

D:=\begin{equation}
\{(x_1,x_2,\cdots,x_{n-1})\ |\ 0<x_l<\alpha_l, l=1,2,\cdots,n-1
\}
\end{equation}

\noindent according to the precision required; (by Theorem 2 and
its proof)

{\bf Step 3.} Take $b:=(b_1,b_2,\cdots,b_{n-1})$ at each gridding
point. Test if $b$ belongs to $\cap _{i=1}^m \Omega ^W_{1a_i}$ by
the following steps:

1) Test if the $(n-1)$-th order polynomial with coefficients
$(1,b_1,b_2,\cdots,b_{n-1})$ belongs to $H^{n-1}$ (by Property 3),

2) For $ i=1,2,\cdots,m,$ test if the polynomial with coefficients
$(a^{(i)}_1-b_1,a^{(i)}_2-b_2,\cdots,
a^{(i)}_{n-1}-b_{n-1},a^{(i)}_n)$ belongs to $H^{n-1}\cup
H^{n-2}$, respectively (by Theorem 2 and its proof),

3) Test if $b$ belongs to $\Omega ^W_{1a_i}, i=1,2,\cdots,m$ .

If all 1), 2), 3) above are satisfied, go to Step 4; otherwise,
move to the next gridding point and test 1), 2), 3) again; (If all
gridding points have been tested, then print "there does not exist
such a $b$ in $\cap _{i=1}^m
 \Omega ^W_{1a_i}$ with the given precision").

{\bf Step 4.} Take a sufficiently small $\varepsilon>0$ such that
$(\varepsilon,1,b_1,b_2,\cdots,b_{n-1})\in \cap _{i=1}^m \Omega
_{a_i}.$ Hence, the $n$-th order polynomial with coefficients
$(\varepsilon,1,b_1,b_2,\cdots,b_{n-1})$ satisfies the design
requirement. (by Theorem 4).

For the low order stable interval polynomial family or low order
stable convex combination, existence of the solution to the design
problem is always guaranteed[15,16].
 Given adequate precision, our method will surely find a polynomial
that satisfies the design requirement. As shown by numerous
examples below, our method is also effective for higher order
polynomial families.

{\bf  Example  1} \ \

Let

\begin{equation}
\begin{array}{ll}
F  =\{ &  a_1(s)=s^4+11s^3+56s^2+88s+1,\\
    &     a_2(s)=s^4+11s^3+56s^2+88s+50,\\
    &     a_3(s)=s^4+89s^3+56s^2+88s+1,\\
    &     a_4(s)=s^4+89s^3+56s^2+88s+50 \}
\end{array}
\end{equation}

\noindent the methods proposed in [8-10,14,15] do not work here.
 Using our method, it is easy to get
$b(s)= s^3+3.3s^2+2.24s+1.76\in \cap _{i=1}^4
 \Omega ^W_{1a_i}$. Then let $c(s):=\varepsilon s^4+b(s)$, where
$\varepsilon>0$ is sufficiently small, e.g., let $\varepsilon \leq
0.3$ (which is determined by Theorem 4),it is easy to check that
the design requirement has been met.

Note that the example above is constructed by overbounding the
line segment in [13] by an interval polynomial family. Thus,
instead of dealing with two vertex polynomials as in [13], we must
now deal with four Kharitonov's vertex polynomials.

In what follows, we will give two more examples of higher order
polynomial families.

{\bf  Example  2} \ \

Let
{\scriptsize{
\begin{eqnarray}
F=\{\hspace{7.5cm}\nonumber\\
a_1(s)=s^7+9s^6+31s^5+71.5s^4+111.5s^3+109s^2+76s+12.5, \hspace{0.2cm}\nonumber\\
a_2(s)=s^7+9.5s^6+31s^5+71s^4+111.5s^3+109.5s^2+76s+12,\hspace{0.2cm}\nonumber\\
a_3(s)=s^7+9s^6+31.5s^5+71.5s^4+111s^3+109s^2+76.5s+12.5,\nonumber\\
a_4(s)=s^7+9.5s^6+31.5s^5+71s^4+111s^3+109.5s^2+76.5s+12 \}
\end{eqnarray}
}}
\noindent It is easy to see that the convex hull $\overline{F}$ of
$F$ is robust stable. Using our method, it is easy to get
$b(s)=s^6+7.2s^5+18.6s^4+42.6s^3+44.4s^2+43.6s+15.2 \in \cap
_{i=1}^4
 \Omega ^W_{1a_i}.$

Then let $c(s):=\varepsilon s^7+b(s)$, where $\varepsilon>0$ is
sufficiently small, e.g., let $\varepsilon \leq 0.1,$ it is easy
to check that the design requirement has been met.

{\bf  Example  3} \ \

Let
{\scriptsize{
\begin{equation}
\begin{array}{ll}
F=\{ &a_1(s)=s^9+11s^8+52s^7+145s^6+266s^5+331s^4\\
& \hspace{3.5cm}+280s^3+155s^2+49s+6,\\
 &    a_2(s)=s^9+11s^8+52s^7+146s^6+265.5s^5+332s^4\\
 & \hspace{3.2cm}+278.5s^3+151s^2+48s+2 \}
\end{array}
\end{equation}
}}

\noindent It can be verified that
$a_2(s)-a_1(s)=s^6-0.5s^5+s^4-1.5s^3-4s^2-s-4$ satisfies the
extended Alternating Hurwitz Minor
Condition$^{[7,21,22]}$£¬namely, it is  a convex direction for
Hurwitz stability$^{[7,21,22]}$. Moreover, it is easy to see that
$a_1(s)$ and $a_2(s)$ are both Hurwitz stable polynomials. Thus,
the convex hull $\overline{F}$ of $F$ is robust
stable$^{[7,21,22]}$. Using our method, it is easy to get
$b(s)=s^8+8.8s^7+41.6s^6+87s^5+159.3s^4+132.4s^3+111.4s^2+30.2s+9.6
\in \Omega ^W_{1a_1}\cap
 \Omega ^W_{1a_2}.$
Thus, let $c(s):=\varepsilon s^9+b(s), \varepsilon>0, \varepsilon
$ sufficiently small, e.g., take $\varepsilon \leq 0.07,$ then the
design requirement has been met.

Note that our design method is also effective when dealing with
discrete time systems. Note also that, in the Examples 1-3, $b(s)$
is not unique. Using our method, we can get all such $b(s)$'s with
given precision.

It should be pointed out that there is hardly any example with
order higher than 4 in the literature. Recently, a sixth-order
example of interval family was given in [23] as follows.
Unfortunately, this example is incorrect.

{\bf  Example  4} \ \

Suppose
{\scriptsize{
\begin{equation}
\begin{array}{ll}
F=\{ &\\ a_1(s)=s^6+0.8s^5+58.06s^4+50.9s^3+1028.5s^2+163.82s+1042.5,\\
      a_2(s)=s^6+1.5s^5+58.06s^4+28.3s^3+1028.5s^2+376.36s+1042.5,\\
      a_3(s)=s^6+0.8s^5+68.62s^4+50.9s^3+755.47s^2+163.82s+3286.7,\\
      a_4(s)=s^6+1.5s^5+68.62s^4+28.3s^3+755.47s^2+376.36s+3286.7 \},
\end{array}
\end{equation}
}}

\noindent find a polynomial $b(s),$ such that
$$p_i(s):=\displaystyle\frac {b(s)}{a_i(s)}\in SPR, i=1,2,3,4.$$

By Definition 1 , Definition 2 and Property 3, a prerequisite of
the robust SPR design problem is that the convex hull
$\overline{F}$ of $F$ is robustly stable. But it is easy to check
that $\overline{F}$ is not robustly stable. In fact, $(1.0446\pm
5.8969i)$ are roots of $a_1(s)$ with positive real part;
$(1.037\pm 4.9835i)$  are roots of $a_2(s)$ with positive real
part; $(0.03291\pm 7.5026i)$ and $(0.68089\pm 2.4933i)$ are roots
of $a_3(s)$ with positive real part;  $(0.87123\pm 2.867i)$ are
roots of $a_4(s)$ with positive real part. Thus, it does not make
sense to consider the robust SPR design in this case.

It should also be pointed out that, for the vertex set of a
general polytopic polynomial family
$F=\{a_i(s)=s^n+\sum_{l=1}^na_l^{(i)}s^{n-l},i=1,2,\cdots,m.\}$,
even if $\overline{F}$ is robustly stable,it is still possible
that there does not exist a polynomial $c(s),$ such that, $
c(s)/f(s)\in $WSPR, for all $f(s)\in \overline{F}$. Namely $\cap
_{i=1}^m
 \Omega ^W_{1a_i}= \phi.$

To see this, let us look at an example of a third order triangle
polynomial family.

{\bf   Proposition 1}$^{[15,16,18]}$\ \

Let $ a(s)=s^3+a_1s^2+a_2s+a_3\in H^3 ,$ then

\begin{center}
$ \begin{array}{ll} \ \ \  \Omega ^W_{1a}=&
\{(x,y)|a_2^2x^2+2(2a_3-a_1a_2)xy+a_1^2y^2\\
&\quad \quad-2a_2a_3x-2a_1a_3y+a_3^2<0\}%
\\
& \cup \{(x,y)|x \leq a_1,a_2x-a_1y-a_3\geq 0, y>0\}. \\
\end{array}$
\end{center}

{\bf  Example  5} \ \

Let

\begin{equation}
\begin{array}{ll}
F=\{ & a_1(s)=s^3+2.6s^2+37s+64,\\
     & a_2(s)=s^3+17s^2+83s+978,\\
     & a_3(s)=s^3+15s^2+28s+415 \}
\end{array}
\end{equation}

\noindent  it is easy to verify that $a_i(s),i=1,2,3,$ are Hurwitz
stable. Moreover, all edges of $\overline{F}$, i.e., $\lambda
a_i(s)+(1-\lambda)a_j(s), \lambda \in [0,1],i,j=1,2,3,$ are also
Hurwitz stable. Therefore, by Edge Theorem$^{[1,6,7,20]}$,
$\overline{F}$ is robustly stable. On the other hand, by a direct
computation
 using Proposition 1,  we have
$\Omega ^W_{1a_1}\cap \Omega ^W_{1a_2}\cap \Omega ^W_{1a_3}=\phi.$
Henceforth,  there does not exist a polynomial $c(s)$
 such that $c(s)/a_i(s)\in $ WSPR,
 $i=1,2,3$  (although $\Omega ^W_{1a_1}\cap \Omega ^W_{1a_2}\neq \phi,
\Omega ^W_{1a_1}\cap \Omega ^W_{1a_3}\neq \phi,$ and $\Omega
^W_{1a_2}\cap \Omega ^W_{1a_3}\neq \phi$).

Note that, in this example, though we have $\Omega ^W_{1a_1}\cap
\Omega ^W_{1a_2}\cap \Omega ^W_{1a_3}=\phi,$
 we do not know whether $\Omega _{a_1}\cap \Omega _{a_2}\cap \Omega _{a_3}$
is an empty set or not. This is a problem deserving further study.

For a fourth ( or higher ) order stable interval polynomial family
( or stable convex combination of two polynomials), does there
exist a polynomial such that their ratios are SPR-invariant? This
is still an open problem$^{[1,8,9,13,14]}$. From our numerous
examples, it seems that such a polynomial can always be found.
Thus, we conjecture that this problem has a positive answer.

\vspace*{1.5\baselineskip} \centerline{\Large References}
\vspace*{1.0\baselineskip}

\def\toto#1#2{\centerline{\hbox to0.4cm{#1\hss}
\parbox[t]{7.6cm}{#2}}\vspace{0.5\baselineskip}}

\toto{[1]} {L. Wang and L. Huang, Extreme Point Results for Strict
Positive Realness of Transfer Function Families, {\it Systems
Science and Mathematical Sciences,} {\bf 7:} 430-439, 1994. }

\toto{[2]} {Y. D. Landau,
 {\it Adaptive Control: The Model Reference Approach.}
Marcel Dekker, New York, 1979.}

\toto{[3]} {V. M. Popov, {\it Hyperstability of Automatic Control
Systems. } Springer-Verlag, New York, 1973.}

\toto{[4]} {L. Wang and L. Huang, Robust Stability of Polynomial
Families and Robust Strict Positive Realness of Rational Function
Families, {\it Int. J. Systems Science,} {\bf 23:} 235-247, 1992.
}

\toto{[5]} {L. Huang,
 {\it Stability Theory,}
Peking University Press, Beijing, 1989.}

\toto{[6]} {S. P. Bhattacharyya, H. Chapellat and L . H . Keel,
{\it Robust Control: The Parametric Approach.} Prentice Hall , New
York, 1995.}

\toto{[7]} {B. R. Barmish, {\it New Tools for Robustness of Linear
Systems.} MacMillan Publishing Company, New York, 1994.}

\toto{[8]}{S. Dasgupta and A. S. Bhagwat, Conditions for designing
strictly positive real transfer functions for adaptive output
error identification. {\it IEEE Trans. on Circuits and Systems,}
{\bf CAS-34:} 731-737, 1987.}

\toto{[9]}{ L. Huang, C. V. Hollot, Robust analysis of strictly
positive real function set. Preprints of {\it The Second
Japan-China Joint Symposium on Systems Control Theory and its
Applications,} 210-220, 1990.}

\toto{[10]}{ C. V. Hollot, L. Huang, Designing strictly positive
real transfer function families: a necessary and sufficient
condition for low degree and structured families. {\it Proc. of
the International Conference on Mathematical Theory of Network and
Systems,} 215-227, 1989.}

\toto{[11]}{H. Chapellat and S. P. Bhattacharyya, On robust
nonlinear stability of interval control systems. {\it IEEE Trans.
on Automatic Control,} {\bf AC-36:} 59-67, 1991.}

\toto{[12]} {L. Wang and L. Huang, Finite Verification of Strict
Positive Realness of Interval Transfer Functions, {\it Science
Bulletin,} {\bf 4:} 262-268, 1991. }

\toto{[13]}{A. Betser and E. Zeheb, Design of robust strictly
positive real transfer functions. {\it IEEE Trans. Circuits and
Systems,} {\bf CAS-40:} 573-580, 1993.}

\toto{[14]}{B. D. O. Anderson, S. Dasgupta, P. Khargonekar, F. J.
Kraus and M. Mansour, Robust strict positive realness:
characterization and construction. {\it IEEE Trans. Circuits and
Systems,} {\bf CAS-37:} 869-876, 1990.}

\toto{[15]} {W. Yu and L. Huang, Necessary and Sufficient
Conditions for Strictly Positive Real Stabilization of Low-Order
Systems, {\it Science Bulletin,} {\bf 43:} 2275-2279, 1998. }

\toto{[16]}
 {W. Yu,
{\it Robust SPR Synthesis and Robust Stability Analysis.} PhD
Thesis, Peking University, Beijing, 1998.}

\toto{[17]} { W. S. Yu and L. Wang, Some Remarks on the Definition
of Strict Positive Realness of Transfer Functions.
 {\it Proceedings of the Conference on Control and Decision}, Nanjing, 1999.}

\toto{[18]}{ W. S. Yu, L. Wang and M. Tan, Complete
characterization of strict positive realness regions in
Coefficient space. {\it Proceedings of IEEE Hong Kong Symposium on
Robotics and Control,} Hong Kong, 1999.}

\toto{[19]} {J. Ackermann, {\it Robust Control: Systems with
Uncertain Physical Parameters.} Springer-Verlag, London, 1993.}

\toto{[20]}{A. C. Bartlett, C. V. Hollot and L. Huang, Root
locations for an entire polytope of polynomial: it suffices to
check the edges. {\it Mathematics of Control, Signals, and
Systems,} {\bf 1:} 61--71, 1988.}

\toto{[21]} {Barmish, B. R. and Kang, I. H., Extreme point results
for robust stability of interval plants: Beyond first order
compensators. {\it Automatica, } 1992, {\bf 28:} 1169-1180.}

\toto{[22]}{Rantzer, A., Stability conditions for polytopes of
polynomials. {\it IEEE Trans. on Automatic Control, } 1992, {\bf
AC-37:} 79-89.}

\toto{[23]}{Marquez, H. J. and Agathoklis, P., On the existence of
robust strictly positive real rational functions. {\it IEEE Trans.
on Circuits and Systems, Part I,} 1998, {\bf CAS-45:} 962-967.}

\end{multicols}
\end{document}